\documentstyle[12pt]{article}
\newcommand{\nc}{\newcommand}
\nc{\al}{\alpha}
\nc{\ald}{{\dot \al}}
\nc{\ba}{\beta_\al}
\nc{\bb}{\beta_\beta}
\nc{\ga}{\g^\al}
\nc{\gb}{\g^\beta}
\nc{\db}{\pa_\beta}
\nc{\dtb}{\delta_\theta^\beta}
\nc{\dab}{{\delta_\al}^\beta}
\nc{\vmab}{V_{-\al}^\beta}
\nc{\vab}{V_\al^\beta}
\nc{\vib}{V_i^\beta}
\nc{\g}{\gamma}
\nc{\G}{\Gamma}
\nc{\D}{\Delta}
\nc{\paj}{P_{-\al}^j}
\nc{\la}{\lambda}
\nc{\La}{\Lambda}
\nc{\var}{\varphi}
\nc{\kvt}{\sqrt{t}}
\nc{\hn}{h^\vee}
\nc{\kn}{k^\vee}
\nc{\pa}{\partial}
\nc{\nn}{\nonumber \\ }
\nc{\hf}{\frac{1}{2}}         
\nc{\dz}{\frac{dz}{2\pi i}}
\nc{\fabc}{{f_{ab}}^c}
\nc{\binomial}[2]{\left (\begin{array}{c} {#1}\\ {#2} \end{array}
\right )}
\nc{\cpp}{{C_+}^+}
\nc{\cnp}{{C_0}^+}
\nc{\cmp}{{C_-}^+}
\nc{\cmn}{{C_-}^0}
\nc{\cmm}{{C_-}^-}
\nc{\vpp}{{V_+}^+}
\nc{\ben}{\begin{equation}}
\nc{\een}{\end{equation}}
\nc{\bea}{\begin{eqnarray}}
\nc{\eea}{\end{eqnarray}}
\nc{\bra}[1]{\langle {#1}|}
\nc{\ket}[1]{|{#1}\rangle}
\newcommand{\Z}{\mbox{$Z\hspace{-2mm}Z$}}
\nc{\C}{\mbox{\hspace{1.24mm}\rule{0.2mm}{2.5mm}\hspace{-2.7mm} C}}
\nc{\Nat}{\mbox{\hspace{.04mm}\rule{0.2mm}{2.8mm}\hspace{-1.5mm} N}}

\nc{\spa}{\hspace{1 cm},\hspace{1 cm}}
\nc{\vs}{\vspace}
\nc{\NP}[1]{Nucl.\ Phys.\ {\bf #1}}
\nc{\PL}[1]{Phys.\ Lett.\ {\bf #1}}
\nc{\CMP}[1]{Commun.\ Math.\ Phys.\ {\bf #1}}
\nc{\PR}[1]{Phys.\ Rev.\ {\bf #1}}
\nc{\PRL}[1]{Phys.\ Rev.\ Lett.\ {\bf #1}}
\nc{\PTP}[1]{Prog.\ Theor.\ Phys.\ {\bf #1}}
\nc{\PTPS}[1]{Prog.\ Theor.\ Phys.\ Suppl.\ {\bf #1}}
\nc{\MPL}[1]{Mod.\ Phys.\ Lett.\ {\bf #1}}
\nc{\IJMP}[1]{Int.\ Jour.\ Mod.\ Phys.\ {\bf #1}}
\nc{\IM}[1]{Invent.\ Math.\ {\bf #1}}
\nc{\SJNP}[1]{Sov. J. Nucl. Phys.\ {\bf #1}}

\begin{document}

\topmargin -5mm
\oddsidemargin 5mm

\begin{titlepage}
\setcounter{page}{0}
\begin{flushright}
November 1997
\end{flushright}

\vs{8mm}
\begin{center}
{\Large Screening Currents in Affine Current Algebra}

\vs{8mm}
{\large J{\o}rgen Rasmussen}\footnote{e-mail address: 
jorgen@celfi.phys.univ-tours.fr}\\[.2cm]
{\em Laboratoire de Math\'{e}mathique et Physique Th\'{e}orique,}\\
{\em Universit\'{e} de Tours,
Parc de Grandmont, F-37200 Tours, France}

\end{center}

\vs{8mm}
\centerline{{\bf{Abstract}}}
\noindent
In this paper screening currents of the second kind are considered.
They are constructed in any affine current algebra for directions 
corresponding to simple roots with multiplicity one in a decomposition of 
the highest root on a set of simple roots. These expressions are precisely 
of the form previously conjectured
to be valid for all directions in general affine current algebras. However, by
working out explicitly the screening currents in the case of $SO(5)$ based
on the Lie algebra $B_2$, it is demonstrated that much more complicated
structures appear in the general case.
In the distinguished representation of affine $OSp(2|2)$ current superalgebra,
the screening current of the second kind in the bosonic direction is also
provided.
\\[.4cm]

\end{titlepage}
\newpage
\renewcommand{\thefootnote}{\arabic{footnote}}
\setcounter{footnote}{0}

\section{Introduction}

Since the work by Wakimoto \cite{Wak} on free field realizations of affine
$SL(2)$ current algebra much effort has been made in obtaining similar
constructions in the general case, a problem in principle solved by Feigin
and Frenkel \cite{FF}. Recently two independent methods have 
led to general and 
explicit solutions \cite{deBF,Ras1,PRY4}. The method used by Petersen, Yu and
the present author \cite{Ras1,PRY4} gives particularly simple and compact
free field realizations and is amenable of generalizations to affine current
superalgebras \cite{Ras2}. 

Free field realizations
enable one in principle to build integral representations for correlators in
conformal field theory \cite{DF,FZ,F,BF,FGPP,PRY1,An}. 
In a recent series of papers 
Petersen, Yu and the present author
have carried out such a study for conformal field theory based on
affine $SL(2)$ current algebra \cite{PRY1}. It turns out that 
screening operators of both the first and the second
kinds are crucial for being able to treat the general case of degenerate 
representations \cite{KK} and in particular
admissible representations \cite{KW}. In that
connection it is also necessary to be able to handle fractional 
powers of free fields. Well defined rules for that have been
established also in \cite{PRY1}. A particular
interest in these techniques is due to their close relationship with 
2D quantum gravity and string theory \cite{HY,AGSY}.

In order to generalize our work on affine $SL(2)$ current algebra to 
higher groups and supergroups, one needs not only free field realizations
of the affine currents but also of the screening currents and 
the primary fields.
These are also provided in \cite{Ras1,PRY4,Ras2}, though in the general
case not for screening currents of the second kind. An expression, originally
written down by Ito \cite{Ito1}, was proposed as a candidate 
for screening currents of the second kind in the case of arbitrary affine
current algebra, and a proof was provided in the case of $SL(r+1)$.

In this paper we shall prove that the proposal is valid in all affine current
algebras for directions
corresponding to simple roots with multiplicity {\em one} in a decomposition
of the highest root on a set of simple roots. This includes all directions 
in $SL(r+1)$, and indeed the proof we present is identical to the one employed
in \cite{Ras1,PRY4} for $SL(r+1)$.
The main result in this paper is the explicit construction
of the screening current of the second kind in the direction $\al_2$
of multiplicity {\em two}
in affine $SO(5)$ current algebra, based on the Lie algebra $B_2$.
The expression clearly demonstrates that much more complicated
structures than previously anticipated appear in the general case. 
The construction is not unique in the form we present it. Rather, it 
involves an infinite summation over a variable which is only restricted to
take on integer-spaced values and may thus be written as $n\in(\Z+a)$, where
$a$ is a free parameter. This does not come as a surprise since our techniques
for handling fractional powers of free ghost fields allow for the introduction 
of adjustable monodromy parameters \cite{PRY1}. The observed non-uniqueness 
is then an indication that our result for the screening current is a 
fractional (and asymptotic) expansion of some expression. Our result
diverges as $n!$ so it can at most be Borel summable.

The screening current of the second kind for the bosonic direction
in the distinguished representation of $OSp(2|2)$ is also provided. This
result does not differ substantially from the result for the screening 
current of the second kind in $OSp(1|2)$ \cite{ERdeS}, hitherto the only 
known screening current of the second kind in affine current superalgebras. 
Generalizations to higher
groups and supergroups are currently being investigated.

The remaining part of this paper is organized as follows. Section 2 serves
to fix notation and reviews some of the main results in \cite{PRY4}. 
In Section 3 we present a proof of the result in directions of simple
roots with multiplicity one.
The expressions for the screening currents of the second kind
in the direction $\al_2$ in affine $SO(5)$ current algebra and the bosonic
direction in affine $OSp(2|2)$ current superalgebra
are then given and proofs are outlined. Section 4 contains concluding remarks.

\section{Notation}
Let {\bf g} be a simple Lie algebra of dim {\bf g} = $d$ and rank {\bf g} = 
$r$.
{\bf h} is a Cartan subalgebra of {\bf g}. The set of (positive) roots
is denoted ($\Delta_+$) $\Delta$, and we write $\al>0$ if $\al\in\Delta_+$.
The simple roots are $\{\al_i\}_{i=1,...,r}$. $\theta$ is the highest root,
while $\al^\vee = 2\al/\al^2$ is the root dual to $\al$. 
Using the triangular decomposition 
\ben
 \mbox{{\bf g}}=\mbox{{\bf g}}_-\oplus\mbox{{\bf h}}\oplus\mbox{{\bf g}}_+
\een
the raising and lowering operators are denoted $e_\al\in$ {\bf g}$_+$ and
$f_\al\in$ {\bf g}$_-$ respectively with $\al\in\Delta_+$, and 
$h_i\in$ {\bf h} are the Cartan operators. 
We let $j_a$ denote an arbitrary Lie algebra element. 
In the Cartan-Weyl basis we have
\ben
 [h_i,e_\al]=(\al_i^\vee,\al)e_\al\spa [h_i,f_\al]=
  -(\al_i^\vee,\al)f_\al
\label{CW}
\een
and
\ben
 \left[e_\al,f_\al\right]=h_\al=G^{ij}(\al_i^\vee,\al^\vee)h_j
\een
where the metric $G_{ij}$ is related to the Cartan matrix $A_{ij}$ as
\ben
 A_{ij}=\al_i^\vee\cdot\al_j=(\al_i^\vee,\al_j)=
  G_{ij}\al_j^2/2
\een
while the Cartan-Killing form $\kappa_{ab}$ is
\ben
   \kappa_{\al,-\beta}=\frac{2}{\al^2}
  \delta_{\al,\beta}\spa \kappa_{ij}=G_{ij}
\een
The Weyl vector $\rho=\hf\sum_{\al>0}\al$ satisfies
$\rho\cdot\al_i^\vee=1$.
We use the convention ${f_{-\al,-\beta}}^{-\g}=-{f_{\al,\beta}}^\g$ so the
standard symmetries of the structure coefficients may be summarized as
\bea
 {f_{-\al,-\beta}}^{-\g}&=&-{f_{\al,\beta}}^\g \spa{f_{\beta,\al}}^\g=
  -{f_{\al,\beta}}^\g\nn
 \frac{{f_{\al,\beta}}^{\al+\beta}}{(\al+\beta)^2}&=&
  \frac{{f_{\beta,-(\al+\beta)}}^{-\al}}{\al^2}=
  \frac{{f_{-(\al+\beta),\al}}^{-\beta}}{\beta^2}
\label{symmf}
\eea
The Dynkin labels $\La_k$ of the weight $\La$ are defined by
\ben
 \La=\La_k\La^{k}\spa \La_k=(\al_k^\vee,\Lambda)
\een
where $\left\{\La^{k}\right\}_{k=1,...,r}$ is the set of fundamental
weights satisfying
\ben
 (\al_i^\vee,\La^{k})=\delta_i^k
\een 
Elements in $\mbox{\bf g}_+$ or vectors in 
representation spaces (see below) are parametrized using ``triangular 
coordinates" denoted by $x^\al$, one for each positive root, thus we 
write general Lie algebra elements in $\mbox{\bf g}_+$ as
\ben
 g_+(x)=x^\al e_\al \in \mbox{\bf g}_+
\een
and the corresponding group elements $G_+(x)$ as
\ben
 G_+(x)=e^{g_+(x)} 
\een
The matrix representation $C(x)$ of $g_+(x)$ in the adjoint 
representation is defined by
\ben
 C_a^b(x)={C(x)_a}^b={(x^\beta C_\beta)_a}^b=-x^\beta {f_{\beta a}}^b
\label{cadj}
\een
and may be block decomposed as
\ben
 C=\left(\begin{array}{lll}\cpp & 0 & 0\\
                 \cnp & 0 & 0\\
                 \cmp & \cmn & \cmm
         \end{array}  \right)
\label{C}
\een
${C_+}^+$ etc are matrices themselves. In ${C_+}^+$ both row and column 
indices are positive roots, in $\cmn$ the row index is a negative root and
the column index is a Cartan algebra index, etc.
Note that $C_\al^\beta(x)$ vanishes unless 
$\al <\beta$, corresponding to ${C_+}^+$ being upper triangular with 
zeros in the diagonal. Similarly, ${C_-}^-$ is lower triangular. 
We will understand ``properly" repeated root indices as in (\ref{cadj})
to be summed over the positive roots.

For the associated affine Lie algebra, the operator product expansion, OPE, 
of the associated currents is
\ben
 J_a(z)J_b(w)=\frac{\kappa_{ab}k}{(z-w)^2}+\frac{\fabc J_c(w)}{z-w}
\label{JaJb}
\een
where regular terms have been omitted. $k$ is the central extension and
$\kn=2k/\theta^2$ is the level.
The Sugawara energy momentum tensor is
\ben
 T(z)=\frac{1}{2t}\kappa^{ab}:J_aJ_b:(z)
\een
where we have introduced the parameter
\ben
 t=\frac{\theta^2}{2}\left(\kn+\hn\right)
\een
and where $\hn$ is the dual Coxeter number. This tensor has central charge
\ben
 c=\frac{\kn d}{\kn+\hn}
\label{c}
\een

The standard free field construction 
\cite{Wak,FF,GMMOS,BMP,KOS,Ito1,Ito,Dot,Ku,ATY}
consists in introducing for every positive 
root $\al>0$, a pair of free bosonic ghost
fields ($\ba,\ga$) of conformal weights (1,0) satisfying the OPE
\ben
 \ba(z)\gb(w)=\frac{\dab}{z-w}
\een
The corresponding energy-momentum tensor is
\ben
 T_{\beta\g}(z)=:\pa\ga(z)\ba(z):
\label{Tbg}
\een
with central charge
\ben
 c_{\beta\g}=d-r
\een
For every Cartan index $i=1,...,r$ one introduces a free scalar boson $\var_i$
with contraction
\ben
 \var_i(z)\var_j(w)=G_{ij}\ln(z-w)
\een
The energy-momentum tensor 
\ben
 T_\var(z)=\hf:\pa\var(z)\cdot\pa\var(z):-\frac{1}{\kvt}\rho\cdot\pa^2\var(z)
\een
has central charge
\ben
 c_\var=r-\frac{\hn d}{\kn+\hn}
\een
This follows from Freudenthal-de Vries strange formula 
$\rho^2=\hn\theta^2d/24$. The combination $T_\var+T_{\beta\g}$ is the free
field realization of the Sugawara energy-momentum tensor.

The vertex operator
\bea
 V_\Lambda(z)&=&:e^{\frac{1}{\kvt}\Lambda\cdot\var(z)}: \nn
\Lambda\cdot\var(z)&=&\Lambda_iG^{ij}\var_j(z)
\eea
has conformal weight
\ben
 \Delta(V_\Lambda)=\frac{1}{2t}(\Lambda,\Lambda+2\rho)
\een
It is also affine primary corresponding to highest weight $\Lambda$. 
In \cite{PRY4},
the explicit general construction is provided
of the full multiplet of primary fields, parametrized by the $x^\al$ 
coordinates. A similar and likewise general construction in the case of
superalgebras is provided in \cite{Ras2}.

\subsection{Differential Operator Realization}

The lowest weight vector in the (dual) representation space $\langle\Lambda|$
is introduced as
\ben
 \langle\Lambda|f_\al=0\spa\langle\Lambda|h_i=\La_i\langle\Lambda|
\een
An arbitrary vector in this representation space is parametrized as
\ben
\bra{\Lambda,x}=\bra{\Lambda}G_+(x)
\een
A differential operator realization $\tilde{J}_a(x,\pa,\Lambda)$ of the simple 
Lie algebra {\bf g} may then be defined by
\ben
\bra{\Lambda,x}j_a=\tilde{J}_a(x,\pa,\Lambda)\bra{\Lambda,x}
\een
with $\pa_\al=\pa_{x^\al}$ denoting partial derivative wrt $x^\alpha$.
It is obvious that these operators
satisfy the Lie algebra commutation relations. From the Gauss 
decomposition of $\langle\Lambda|G_+(x)e^{tj_a}$, for $t$ small,
\bea
 \langle\Lambda|G_+(x)\exp(te_\al)&=&\langle\Lambda|\exp\left( 
  x^\g e_\g+t\vab(x)e_\beta+{\cal O}(t^2)\right)\nn
 &=&\langle\Lambda|\exp\left(t\vab(x)\db+{\cal O}(t^2)\right)G_+(x)\nn
 \langle\Lambda|G_+(x)\exp(th_i)&=&\langle\Lambda|\exp\left(th_i\right)
  \exp\left(x^\g e_\g+tV_i^\beta(x)e_\beta+{\cal O}(t^2)\right)\nn
 &=&\langle\Lambda|\exp\left(t\left(V_i^\beta(x)\db+\Lambda_i\right)
  +{\cal O}(t^2)\right)G_+(x)\nn 
 \langle\Lambda|G_+(x)\exp(tf_\al)&=&\langle\Lambda|\exp\left(
  tQ_{-\al}^{-\beta}(x)f_\beta+{\cal O}(t^2)\right)\exp\left(
  t\paj(x)h_j+{\cal O}(t^2)\right)\nn
 &\cdot&\exp\left(x^\g e_\g+t\vmab(x)e_\beta+
  {\cal O}(t^2)\right)\nn
 &=&\langle\Lambda|
  \exp\left(t\left(\paj(x)\Lambda_j+\vmab(x)\db\right)
  +{\cal O}(t^2)\right)G_+(x)
\label{Gauss1}
\eea
it follows that the differential operator realization is of the form
\bea
\tilde{E}_\al(x,\pa)&=&\vab(x)\db\nn
\tilde{H}_i(x,\pa,\Lambda)&=&\vib(x)\db+\Lambda_i\nn
\tilde{F}_\al(x,\pa,\Lambda)&=&\vmab(x)\db+\paj(x)\Lambda_j
\label{defVP}
\eea
and one finds \cite{Ras1,PRY4}
\bea
 \vab(x)&=&\left[B(C(x))\right]_\al^\beta\nn
 \vib(x)&=&-\left[C(x)\right]_i^\beta \nn
 \vmab(x)&=&\left[e^{-C(x)}\right]_{-\al}^\g\left[B(-C(x))\right]_\g^\beta\nn
 \paj(x)&=&\left[e^{-C(x)}\right]_{-\al}^j \nn
  Q_{-\al}^{-\beta}(x)&=&\left[e^{-C(x)}\right]_{-\al}^{-\beta}
\label{VPQ}
\eea
$B$ is the generating function for the Bernoulli numbers
\bea
  B(u)&=&\frac{u}{e^u-1}=\sum_{n\geq 0}\frac{B_n}{n!}u^n\nn
  B^{-1}(u)&=&\frac{e^u-1}{u}=\sum_{n \geq 0}\frac{1}{(n+1)!}u^n
\label{Ber}
\eea
The matrix functions (\ref{VPQ}) are defined in terms of universal
power series expansions, valid for any Lie algebra, but ones that truncate 
giving rise to finite polynomials of which the explicit forms depend on the
Lie algebra in question.

The differential screening operators are defined by
\bea
 \exp\{-te_\al\}G_+(x)&=&\exp\{tS_\al(x,\pa)+{\cal O}(t^2)\}G_+(x)\nn
 S_\al(x,\pa)&=&S_\al^\beta(x)\db
\eea
and are seen to satisfy
\ben
 S_\al(x,\pa)=\tilde{E}_\al(-x,-\pa)
\een
so that
\ben
 S_\al^\beta(x)=-\left[B(-C(x))\right]_\al^\beta
\label{SB}
\een
The screening currents of the first kind
are constructed using these polynomials for simple roots as building blocks,
see below.

\subsection{Free Field Realization}

The free field realization is obtained from the differential 
operator realization $\left\{\tilde{J}_a\right\}$ by the substitution
\cite{FF,GMMOS,BMP,Ito,Ku,ATY}
\ben
 \pa_\al\rightarrow\beta_\al(z)\spa x^\al\rightarrow\g^\al(z)
  \spa\Lambda_i\rightarrow\kvt\pa\varphi_i(z)
\een
and a subsequent addition of a 
normal ordering contribution or anomalous term to the lowering part,
giving rise to the following form of the free field realization
\bea
 E_\al(z)&=&:\vab(\g(z))\beta_\beta(z):\nn
 H_i(z)&=&:\vib(\g(z))\beta_\beta(z):+\kvt\pa\varphi_i(z)\nn
 F_\al(z)&=&:\vmab(\g(z))\beta_\beta(z):+\kvt\pa\varphi_j(z)\paj(\g(z))
  +\pa\g^\beta(z)F_{\al\beta}(\g(z))\nn
 \Delta(J_a)&=&1
\label{Wakimoto}
\eea
In \cite{PRY4} (see also \cite{Ras1,deBF} and \cite{Ras2}) 
the explicit form of $F_{\al\beta}$ is found to be
\bea
 F_{\al\beta}(\g)&=&\frac{2k}{\al^2}
  \left((V_+^+(\g))^{-1}\right)_\beta^\al
  +\left((V_+^+(\g))^{-1}\right)_\beta^\mu\pa_\sigma V_\mu^\g(\g)\pa_\g 
  V_{-\al}^\sigma(\g)
\label{anomal}
\eea
where
\bea
 (V_+^+(\g))^{-1}&=&B(C_+^+(\g))^{-1}\nn
 &=&\sum_{n\geq0}\frac{1}{(n+1)!}(C_+^+(\g))^n
\eea

\subsection{Screening Currents of the First Kind}

A screening current is conformally primary of weight 1 and has the property
that the singular part of the OPE with an affine current is a total 
derivative. These properties ensure that integrated
screening currents (screening charges) may be inserted into correlators
without altering the conformal or affine Ward identities. This in turn makes 
them very useful in construction of correlators, see e.g. 
\cite{DF,BF,Dot,PRY1,Ras1}.
The best known screening currents \cite{FF,BMP,Ito1,Ku,ATY,deBF,Ras1,PRY4} 
are the following denoted screening currents of the first kind, one for each 
simple root
\bea
 s_j(w)&=&:S_{\al_j}^\al
  (\g(w))\beta_\al(w)::e^{-\frac{1}{\kvt}\al_j\cdot\varphi(w)}: \nn
 \al_j\cdot\var(w)&=&\frac{\al_j^2}{2}\var_j(w)
\label{sj}
\eea
satisfying
\bea
  E_\al(z)s_j(w)&=&0\nn
  H_i(z)s_j(w)&=&0\nn
  F_\al(z)s_j(w)&=&\frac{\pa}{\pa w}\left(
    \frac{-2t/\al_j^2}{z-w}Q_{-\al}^{-\al_j}(\g(w))
  :e^{-\frac{1}{\kvt}\al_j\cdot\var(w)}:\right)\nn
  T(z)s_j(w)&=&\frac{\pa}{\pa w}\left(\frac{1}{z-w}s_j(w)\right)
\eea
We shall use the terminology that $s_j(z)$ is the screening current
of the first kind in the direction $\al_j$.

\section{Screening Currents of the Second Kind}

In \cite{BO} Bershadsky and Ooguri found a second screening current in the case
of $SL(2)$
\ben
 \tilde{s}(w)=(-\beta(w))^{-(k^\vee+h^\vee)}:e^{\kvt\var_1(w)}:
\een
Since it involves a generically non-integer power of the free ghost field
$\beta$, discussions on its interpretation remained only partly
successful. However, in \cite{PRY1,Ras1} by Petersen, Yu and the present
author it is demonstrated how techniques of fractional calculus 
provide a solution. As a result we were able to render the free field 
realization applicable of producing integral representations of $N$-point
chiral blocks for degenerate representations and in particular for
admissible representations.

The problem of extending the construction of Bershadsky and Ooguri
to higher groups is discussed in \cite{Ras1,PRY4} where a general proposal
is studied, originally written down by Ito \cite{Ito1}. In the case of
$SL(r+1)$ we there presented a proof of this proposal but for general
groups it has remained a conjecture. However, Proposition 1
below demonstrates that in general
the validity of that proposal is restricted to certain
directions. We shall use the following notation:
\ben
 \theta=\sum_{i=1}^r a^i\al_i
\een
is the decomposition of the highest root $\theta$ on the space of simple
roots, while
\bea
 \tilde{s}_j(w)&=&\tilde{S}_j(w):e^{\kvt\var_j(w)}:\nn
 J_a(z)\tilde{s}_j(w)&=&\frac{\pa}{\pa w}
  \left(\frac{1}{z-w}R_{a,\al_j}(w):e^{\kvt\var_j(w)}:\right) 
\label{defR}
\eea
is the contraction between the affine current $J_a(z)$ and the screening 
current of the second kind $\tilde{s}_j(w)$ in the direction 
$\al_j$. $\tilde{S}_j(w)$ and $R_{a,\al_j}(w)$ are assumed to be 
functionals of the ghost fields and derivatives thereof.\\[.2cm]
{\bf Proposition 1}\\[.2cm]
In the direction of the simple root $\al_j$ of multiplicity one, $a^j=1$,
the screening current of the second kind is given by
\ben
 \tilde{s}_j(w)=:\left(S_{\al_j}^\sigma(\g(w))\beta_\sigma(w)\right)^{-
  2t/\al_j^2}::e^{\kvt\var_j(w)}:
\label{screenmultone}
\een
and it produces
\bea
 R_{\al,\al_j}&=&R_{i,\al_j}=0\nn
 R_{-\al,\al_j}&=&-\frac{2t}{\al_j^2}:Q_{-\al}^{-\al_j}(\g)\left(
  S_{\al_j}^\sigma(\g)\beta_\sigma\right)^{-2t/\al_j^2-1}:\nn
 \D(\tilde{s}_j)&=&1
\label{Rs}
\eea
{\bf Proof}\\[.2cm]
First one computes all possible (multiple) contractions between a generator
($J_a$ or $T$) and the screening current $\tilde{s}_j$. Comparisons with
(\ref{Rs}) will then yield
a set of relations among the polynomials $V$, $P$, $Q$ and $S$. These
are then proven using the various classical and quantum polynomial 
identities which follow from the fact that indeed (\ref{defVP}) and 
(\ref{Wakimoto}) constitute a differential operator realization and a free 
field realization, respectively. \cite{Ras1,PRY4,Ras2} may be consulted
for details and lists of polynomial identities. In this final part of the
proof, the essential point is that for $a^j=1$ certain "contracted
sequences" of polynomials like
$S_{\al_j}^\g\pa_\g S_{\al_j}^\sigma$ and 
$S_{\al_j}^{\g_1}S_{\al_j}^{\g_2}S_{\al_j}^{\g_3}\pa_{\g_1}\pa_{\g_2}
\pa_{\g_3}V_{-\al}^\beta$ all vanish; that is when $\al_j$ net appears more 
than once in the lower indices. This simple rule follows immediately from 
the expressions (\ref{VPQ}), (\ref{SB}) 
and (\ref{anomal}) and the definition of the matrix $C$
(\ref{cadj}) and (\ref{C}). In general $S_{\al_j}^\g\pa_\g 
S_{\al_j}^\sigma\pa_\sigma V_{-\al}^\beta$ will not vanish. It is similarly
obvious that contracted sequences vanish if the sum of roots in the upper
indices is less than the net sum of roots in the lower indices.\\
$\Box$\\
Note that the 3rd order pole in the OPE $T(z)\tilde{s}_j(w)$ is proportional
to $S_{\al_j}^\sigma\pa_\sigma S_{\al_j}^\beta$ which is generically
non-vanishing for $a^j>1$, showing the limited validity of the expression
(\ref{screenmultone}).

\subsection{Case of $SO(5)$}

Here we shall consider the case $SO(5)$
where the affine current algebra is based on the
simple Lie algebra $B_2$. This Lie algebra has rank $r=2$, 4 positive roots and
dimension 10, while the dual Coxeter number is $h^\vee=3$. We shall use the
notation where $\al_{11}=\al_1+\al_2$ and
\ben
 \theta=\al_1+2\al_2
\label{thetaB2}
\een
denote the two non-simple positive roots. The normalization
of the root system is
such that the highest root $\theta$ has length squared $\theta^2=\al_1^2
=2\al_2^2=2\al_{11}^2=2$.
The Cartan matrix and the Cartan-Killing form are given by
\bea
 A_{11}=2\spa A_{12}=-1&\spa&A_{21}=-2\spa A_{22}=2\nn
 G_{11}=2\spa G_{12}=-2&\spa&G_{21}=-2\spa G_{22}=4\nn
 &\kappa_{\al,-\al}=2/\al^2&
\eea
The remaining structure coefficients are (up to the symmetries (\ref{CW}) and
(\ref{symmf}))
\bea
 {f_{\al_1,\al_2}}^{\al_{11}}=1&\spa&{f_{\al_2,\al_{11}}}^{\theta}=2\nn
 {f_{1,\al_{11}}}^{\al_{11}}=1&\spa&{f_{2,\al_{11}}}^{\al_{11}}=0\nn
 {f_{1,\theta}}^\theta=0&\spa&{f_{2,\theta}}^{\theta}=2\nn
 {f_{\al_{11},-\al_{11}}}^1=2&\spa&{f_{\al_{11},-\al_{11}}}^2=1\nn
 {f_{\theta,-\theta}}^1=1&\spa&{f_{\theta,-\theta}}^2=1
\eea
The differential operator realization is worked out to be
\bea
 e_{\al_1}(x)&=&\pa_1-\frac{1}{2}x^2\pa_{11}-\frac{1}{6}x^2x^2\pa_\theta\nn
 e_{\al_2}(x)&=&\pa_2+\frac{1}{2}x^1\pa_{11}+
  \left(\frac{1}{6}x^1x^2-x^{11}\right)\pa_\theta\nn
 e_{\al_{11}}(x)&=&\pa_{11}+x^2\pa_\theta\nn
 e_\theta(x)&=&\pa_\theta\nn
 h_1(x)&=&-2x^1\pa_1+x^2\pa_2-x^{11}\pa_{11}+\La_1\nn
 h_2(x)&=&2x^1\pa_1-2x^2\pa_2-2x^\theta\pa_\theta+\La_2\nn
 f_{\al_1}(x)&=&-x^1x^1\pa_1+\left(\hf x^1x^2-x^{11}\right)\pa_2
  -\hf x^1\left(\hf x^1x^2+x^{11}\right)\pa_{11}\nn
 &+&\frac{1}{3}x^1x^2x^{11}\pa_\theta+x^1\La_1\nn
 f_{\al_2}(x)&=&2\left(\hf x^1x^2+x^{11}\right)\pa_1-x^2x^2\pa_2
  +\left(\frac{1}{3}x^1x^2x^2-x^\theta\right)\pa_{11}\nn
 &-&x^2\left(\frac{1}{3}x^2x^{11}+
    x^\theta\right)\pa_\theta+x^2\La_2\nn
 f_{\al_{11}}(x)&=&-2x^1\left(\hf x^1x^2+x^{11}\right)\pa_1+\left(\frac{2}{3}
  x^1x^2x^2-x^2x^{11}+x^\theta\right)\pa_2\nn
 &+&\left(-\frac{5}{12}x^1x^1x^2x^2
  -\hf x^1x^2x^{11}+\hf x^1x^\theta-x^{11}x^{11}\right)\pa_{11}\nn
 &+&\left(\frac{1}{36}x^1x^1x^2x^2x^2+\hf x^1x^2x^2x^{11}+\frac{1}{6}
   x^1x^2x^\theta-x^{11}x^\theta\right)\pa_\theta\nn
 &+&2\left(\hf x^1x^2+x^{11}
     \right)\La_1-\left(\hf x^1x^2-x^{11}\right)\La_2\nn
 f_\theta(x)&=&\left(\frac{1}{4}x^1x^1x^2x^2+x^1x^2x^{11}+x^{11}x^{11}\right)
  \pa_1-x^2\left(\frac{1}{6}x^1x^2x^2+x^\theta\right)\pa_2\nn
 &+&\left(\frac{1}{8}x^1x^1x^2x^2x^2+\frac{1}{3}x^1x^2x^2x^{11}+\hf
  x^2x^{11}x^{11}-x^{11}x^\theta\right)\pa_{11}\nn
 &-&\left(\frac{1}{72}x^1x^1x^2x^2x^2x^2+\frac{1}{6}x^1x^2x^2x^2x^{11}+
  \frac{1}{6}x^2x^2x^{11}x^{11}+x^\theta x^\theta\right)\pa_\theta\nn
 &+&\left(-\frac{1}{3}x^1x^2x^2-x^2x^{11}+x^\theta\right)\La_1
  +\left(\frac{1}{6}x^1x^2x^2+x^\theta\right)\La_2
\eea
Here we have introduced the simplifying notation $x^{11}=x^{\al_{11}},\
\pa_1=\pa_{\al_1}$ etc. 
In the following we shall also need the analogous
abbreviations $\g^2(z)=\g^{\al_2}(z),\ \beta_{11}(z)=\beta_{\al_{11}}(z)$ etc. 
Furthermore, the differential screening operators are
\bea
 S_{\al_1}^\beta(x)&=&-\pa_1-\hf x^2\pa_{11}+\frac{1}{6}x^2x^2\pa_\theta\nn
 S_{\al_2}^\beta(x)&=&-\pa_2+\hf x^1\pa_{11}-\left(\frac{1}{6}x^1x^2
  +x^{11}\right)\pa_\theta
\eea
The (generalized) Wakimoto free field realization of the associated affine
current algebra becomes
\bea
 E_{\al_1}&=&\beta_1-\hf\g^2\beta_{11}-\frac{1}{6}\g^2\g^2\beta_\theta\nn
 E_{\al_2}&=&\beta_2+\hf\g^1\beta_{11}+\left(\frac{1}{6}\g^1\g^2-\g^{11}\right)
  \beta_\theta\nn
 E_{\al_{11}}&=&\beta_{11}+\g^2\beta_\theta\nn
 E_\theta&=&\beta_\theta\nn
 H_1&=&-2:\g^1\beta_1:+:\g^2\beta_2:-:\g^{11}\beta_{11}:+\kvt\pa\var_1\nn
 H_2&=&2:\g^1\beta_1:-2:\g^2\beta_2:-2:\g^\theta\beta_\theta:+\kvt\pa\var_2\nn
 F_{\al_1}&=&-:\g^1\g^1\beta_1:+:\left(\hf\g^1\g^2-\g^{11}\right)\beta_2:
  -\hf:\g^1\left(\hf \g^1\g^2+\g^{11}\right)\beta_{11}:\nn
 &+&\frac{1}{3}\g^1\g^2\g^{11}\beta_\theta
  +\kvt\g^1\pa\var_1+\left(k+\frac{1}{2}\right)\pa\g^1\nn
 F_{\al_2}&=&2:\left(\hf\g^1\g^2+\g^{11}\right)\beta_1:-:\g^2\g^2\beta_2:
  +\left(\frac{1}{3}
  \g^1\g^2\g^2-\g^\theta\right)\beta_{11}\nn
 &-&:\g^2\left(\frac{1}{3}\g^2\g^{11}+\g^\theta\right)\beta_\theta:
  +\kvt\g^2\pa\var_2+2(k+1)\pa\g^2\nn
 F_{\al_{11}}&=&-2:\g^1\left(\hf\g^1\g^2+\g^{11}\right)\beta_1:+
   :\left(\frac{2}{3}\g^1\g^2\g^2-\g^2\g^{11}+\g^\theta\right)\beta_2:\nn
 &+&:\left(-\frac{5}{12}\g^1\g^1\g^2\g^2-\hf\g^1\g^2\g^{11}+\hf\g^1\g^\theta
  -\g^{11}\g^{11}\right)\beta_{11}:\nn
 &+&:\left(\frac{1}{36}\g^1\g^1\g^2\g^2\g^2+\hf\g^1\g^2\g^2\g^{11}+\frac{1}{6}
  \g^1\g^2\g^\theta-\g^{11}\g^\theta\right)\beta_\theta:\nn
 &+&2\kvt\left(\hf\g^1\g^2+\g^{11}\right)\pa\var_1-\kvt\left(\hf\g^1\g^2
  -\g^{11}\right)\pa\var_2\nn
 &+&\left(k+\frac{2}{3}\right)\g^2\pa\g^1-\left(k+\frac{11}{6}\right)
  \g^1\pa\g^2+(2k+1)\pa\g^{11}\nn
 F_{\theta}&=&:\left(\frac{1}{4}\g^1\g^1\g^2\g^2+\g^1\g^2\g^{11}+\g^{11}\g^{11}
  \right)\beta_1:-:\g^2\left(\frac{1}{6}\g^1\g^2\g^2+\g^\theta\right)\beta_2:\nn
 &+&:\left(\frac{1}{8}\g^1\g^1\g^2\g^2\g^2+\frac{1}{3}\g^1\g^2\g^2\g^{11}
  +\hf\g^2\g^{11}\g^{11}-\g^{11}\g^\theta\right)\beta_{11}:\nn
 &-&:\left(\frac{1}{72}\g^1\g^1\g^2\g^2\g^2\g^2+\frac{1}{6}\g^1\g^2\g^2\g^2
  \g^{11}+\frac{1}{6}\g^2\g^2\g^{11}\g^{11}+\g^\theta\g^\theta\right)
  \beta_\theta:\nn
 &+&\kvt\left(-\frac{1}{3}\g^1\g^2\g^2-\g^2\g^{11}+\g^\theta\right)\pa\var_1
  +\kvt\left(\frac{1}{6}\g^1\g^2\g^2+\g^\theta\right)\pa\var_2\nn
 &-&\frac{1}{3}\left(k+\hf\right)\g^2\g^2\pa\g^1+\frac{1}{3}\left(k+\frac{5}{2}
  \right)\g^1\g^2\pa\g^2+(k+2)\g^{11}\pa\g^2\nn
 &-&(k+1)\g^2\pa\g^{11}+k\pa\g^\theta
\eea
For notational reasons we have left out the arguments 
which are the same for all the fields.
The screening currents of the first kind are found to be
\bea
 s_1(w)&=&\left(-\beta_1(w)-\hf\g^2(w)\beta_{11}(w)+\frac{1}{6}\g^2(w)
  \g^2(w)\beta_\theta(w)\right):e^{-\var_1(w)/\kvt}:\nn
 s_2(w)&=&\left(-\beta_2(w)+\hf\g^1(w)\beta_{11}(w)-\left(
  \frac{1}{6}\g^1(w)\g^2(w)+\g^{11}(w)\right)\beta_\theta(w)\right)\nn
 &\cdot&:e^{-\var_2(w)/(2\kvt)}:
\eea

Due to the decomposition (\ref{thetaB2}), it follows from
Proposition 1 that
the screening current $\tilde{s}_1(w)$ is given by 
\bea
 \tilde{s}_1(w)&=&\left(-\beta_1(w)-\hf\g^2(w)\beta_{11}(w)
  +\frac{1}{6}\g^2(w)\g^2(w)\beta_\theta(w)\right)^{-t}:e^{\kvt\var_1(w)}:
\eea 
Note that in this case the normal ordering of the $\beta\g$ part is not
necessary. It may be checked explicitly that the OPE's with the
generators $\{J_a(z)\}$ and $T(z)$ produce
\bea
 R_{-\al_1,\al_1}&=&-t\left(-\beta_1-\hf\g^2\beta_{11}
  +\frac{1}{6}\g^2\g^2\beta_\theta\right)^{-t-1}\nn
 R_{-\al_{11},\al_1}&=&-2t\g^2\left(-\beta_1-\hf\g^2\beta_{11}
  +\frac{1}{6}\g^2\g^2\beta_\theta\right)^{-t-1}\nn
 R_{-\theta,\al_1}&=&t\g^2\g^2
  \left(-\beta_1-\hf\g^2\beta_{11}
  +\frac{1}{6}\g^2\g^2\beta_\theta\right)^{-t-1}\nn 
 R_{-\al_2,\al_1}&=&R_{\al,\al_1}=R_{i,\al_1}=0\nn
 \D(\tilde{s}_1)&=&1
\eea
These expressions comply with the general statement (\ref{Rs}) since 
\bea
 Q_{-\al_1}^{-\al_1}(x)=1&\spa&Q_{-\al_{11}}^{-\al_1}(x)=2x^2\nn
 Q_{-\al_2}^{-\al_1}(x)=0&\spa&Q_{-\theta}^{-\al_1}(x)=-x^2x^2
\eea
{\bf Proposition 2}\\[.2cm]
The screening current $\tilde{s}_2(w)$ is given by
\bea
 \tilde{s}_2(w)&=&\sum_{n}C_n:\left(-\frac{1}{3}\left(2\pa\g^1(w)
   \beta_\theta(w)-\g^1(w)\pa\beta_\theta(w)\right)\right)^n\nn
  &\cdot&\left(-\beta_2(w)
   +\hf\g^1(w)\beta_{11}(w)-\left(\g^{11}(w)
   +\frac{1}{6}\g^1(w)\g^2(w)\right)\beta_\theta(w)
   \right)^{-2t-2n}:\nn
  &\cdot&:e^{\kvt\var_2(w)}:\nn
  C_n&=&\frac{1}{2^nn!}\frac{(-2t)!}{(-2t-2n)!}
\label{sB2}
\eea
and produces
\bea
 R_{-\al_2,\al_2}&=&\sum_{n}C_n(-2t-2n):\left(-\frac{1}{3}\left(2\pa\g^1
   \beta_\theta-\g^1\pa\beta_\theta\right)\right)^n\nn
  &\cdot&\left(-\beta_2+\hf\g^1\beta_{11}-\left(\g^{11}
   +\frac{1}{6}\g^1\g^2\right)\beta_\theta\right)^{-2t-2n-1}:\nn
 R_{-\al_{11},\al_2}&=&\sum_{n}C_n(2t+2n):\g^1\left(-\frac{1}{3}\left(2\pa\g^1
   \beta_\theta-\g^1\pa\beta_\theta\right)\right)^n\nn
  &\cdot&\left(-\beta_2+\hf\g^1\beta_{11}-\left(\g^{11}
   +\frac{1}{6}\g^1\g^2\right)\beta_\theta\right)^{-2t-2n-1}:\nn
 R_{-\theta,\al_2}&=&\sum_{n}C_nn:\pa\g^1 
  \left(-\frac{1}{3}\left(2\pa\g^1
   \beta_\theta-\g^1\pa\beta_\theta\right)\right)^{n-1}\nn
  &\cdot&\left(-\beta_2+\hf\g^1\beta_{11}-\left(\g^{11}
   +\frac{1}{6}\g^1\g^2\right)\beta_\theta\right)^{-2t-2n}:\nn
  &+&\sum_{n}C_n(-2t-2n):\left(\hf\g^1\g^2+\g^{11}\right)
   \left(-\frac{1}{3}\left(2\pa\g^1
   \beta_\theta-\g^1\pa\beta_\theta\right)\right)^n\nn
  &\cdot&\left(-\beta_2+\hf\g^1\beta_{11}-\left(\g^{11}
   +\frac{1}{6}\g^1\g^2\right)\beta_\theta\right)^{-2t-2n-1}:\nn
 R_{-\al_1,\al_2}&=&R_{\al,\al_2}=R_{i,\al_2}=0\nn
 \D(\tilde{s}_2)&=&1
\label{RsB2}
\eea
The summation over $n$ is not restricted to be a summation over integers,
but is given by $n\in(\Z+a)$ where $a$ is a free parameter. For $b$
non-integer, the factorial $b!$ is defined by the gamma-function
$b!=\G(b+1)$.\\[.2cm]
{\bf Proof}\\[.2cm]
The strategy is straightforward, though very tedious. It amounts to compute
all possible (multiple) contractions between a generator and the screening 
current for all generators and then reduce the expressions to obtain the 
total derivatives given above, (\ref{defR}) and (\ref{RsB2}). 
In the process we employ the recursion relation
\ben
 (-2t-2n)(-2t-2n-1)C_n=2(n+1)C_{n+1}
\een
Note that this relation is valid also for $n$ non-integer; that is for
$a$ non-integer.\\
$\Box$\\
The non-uniqueness of the expression (\ref{sB2})
parameterized by $a$ is believed to
be an artifact of our representation of the result. In \cite{PRY1,Ras1} 
(see also \cite{Sch}) it is discussed how fractional (and asymptotic)
expansions like
\ben
 (1+z)^{u}_{(a)}=\sum_{n\in\Z}\binomial{u}{n+a}z^{n+a}
\een
are relevant when computing contractions between ghost fields raised to 
non-integer powers. In that connection $a$ plays the role of an
adjustable monodromy parameter. We believe that a similar situation is 
encountered here and that our result for the screening current is merely
a fractional expansion of some expression. Since our result diverges as
$n!$ it can at most be an asymptotic (fractional) expansion,
or a Borel summable expression.

\subsection{Case of $OSp(2|2)$}

Here we shall present the screening current of the second kind in the bosonic 
direction in the distinguished representation of affine $OSp(2|2)$ current
superalgebra. The notation and results used here are taken from
\cite{Ras2} which may be consulted for further details and 
general results in affine current superalgebras.
 
The Lie superalgebra $A(1,0)$, which is isomorphic to $osp(2|2)$, $C(2)$ and
$sl(1|2)\simeq sl(2|1)$, has rank $r=2$ and 3 positive roots, while
the dual Coxeter number is $\hn=1$ and dim$(\mbox{{\bf g}}_0)=$ 
dim$(\mbox{{\bf g}}_1)=4$.
Here we choose the distinguished set of simple roots consisting of one even
simple root $\al_1$ and one odd simple root $\ald_2$. The remaining and
non-simple root $\ald=\al_1+\ald_2$ is then odd. 
The Weyl vector is $\rho=-\ald_2$, while the
non-vanishing elements of the Cartan-Killing form are 
\bea
 G_{11}=2\spa G_{12}&=&G_{21}=-1\spa G_{22}=0\nn
 \kappa_{\al_1,-\al_1}&=&\kappa_{\ald_2,-\ald_2}=\kappa_{\ald,-\ald}=1
\label{A10G}
\eea
such that the Cartan matrix is given by $A_{ij}=G_{ij}$. 
The remaining non-vanishing structure coefficients are found to be
\bea
 {f_{\al_1,-\al_1}}^1=1&\spa&{f_{\ald_2,-\ald_2}}^2=1\nn
 \ald(H_1)={f_{1,\ald}}^\ald=1&\spa&\ald(H_2)={f_{2,\ald}}^\ald=-1\nn
 {f_{\ald,-\ald}}^1={f_{\ald,-\ald}}^2=1
  &\spa&{f_{\pm\ald_2,\mp\ald}}^{\mp\al_1}=1\nn
 {f_{\pm\al_1,\pm\ald_2}}^{\pm\ald}=\pm1&\spa&
  {f_{\pm\al_1,\mp\ald}}^{\mp\ald_2}=\mp1
\label{A10f}
\eea
The differential operator realization $\{\tilde{J}_A\}$ is 
\bea
 \tilde{H}_1(x,\theta,\pa,\La)&=&-2x^{\al_1}\pa_{\al_1}+\theta^{\ald_2}
  \pa_{\ald_2}-\theta^\ald\pa_\ald+\La_1\nn
 \tilde{H}_2(x,\theta,\pa,\La)&=&x^{\al_1}\pa_{\al_1}+\theta^\ald\pa_\ald
  +\La_2\nn
 \tilde{E}_{\al_1}(x,\theta,\pa)&=&\pa_{\al_1}-\hf\theta^{\ald_2}\pa_\ald\nn
 \tilde{F}_{\al_1}(x,\theta,\pa,\La)&=&-x^{\al_1}x^{\al_1}\pa_{\al_1}
  +\left(\hf x^{\al_1}\theta^{\ald_2}-\theta^\ald\right)\pa_{\ald_2}\nn
 &-&\hf x^{\al_1}\left(\hf x^{\al_1}\theta^{\ald_2}+\theta^\ald\right)\pa_\ald
  +x^{\al_1}\La_1\nn
 \tilde{e}_{\ald_2}(x,\theta,\pa)&=&\pa_{\ald_2}+\hf x^{\al_1}\pa_\ald\nn
 \tilde{f}_{\ald_2}(x,\theta,\pa,\La)&=&\left(\hf x^{\al_1}\theta^{\ald_2}
  +\theta^\ald\right)\pa_{\al_1}+\hf\theta^{\ald_2}\theta^\ald\pa_\ald
  +\theta^{\ald_2}\La_2\nn
 \tilde{e}_{\ald}(x,\theta,\pa)&=&\pa_\ald\nn
 \tilde{f}_{\ald}(x,\theta,\pa,\La)&=&-x^{\al_1}\left(\hf x^{\al_1}
  \theta^{\ald_2}+\theta^\ald\right)\pa_{\al_1}-\theta^{\ald_2}\theta^\ald
  \pa_{\ald_2}\nn
 &+&\left(\hf x^{\al_1}\theta^{\ald_2}+\theta^\ald\right)\La_1-\left(
  \hf x^{\al_1}\theta^{\ald_2}-\theta^\ald\right)\La_2
\eea
The (generalized) Wakimoto free field realization of the associated 
affine current superalgebra becomes
\bea
 H_1(z)&=&-2:\g(z)\beta(z):+:c(z)b(z):-:C(z)B(z):+\kvt\pa\var_1(z)\nn
 H_2(z)&=&:\g(z)\beta(z):+:C(z)B(z):+\kvt\pa\var_2(z)\nn
 E_{\al_1}(z)&=&\beta(z)-\hf c(z)B(z)\nn
 F_{\al_1}(z)&=&-:\g^2(z)\beta(z):+:\left(\hf\g(z)c(z)-C(z)\right)b(z):\nn
 &-&\hf:\g(z)\left(\hf\g(z)c(z)+C(z)\right)B(z):
  +\kvt\g(z)\pa\var_1(z)+(k-1/2)\pa\g(z)\nn
 e_{\ald_2}(z)&=&b(z)+\hf\g(z)B(z)\nn
 f_{\ald_2}(z)&=&:\left(\hf\g(z)c(z)+C(z)\right)\beta(z):+\hf c(z):C(z)B(z):\nn
 &+&\kvt c(z)\pa\var_2(z)+(k+1/2)\pa c(z)\nn
 e_\ald(z)&=&B(z)\nn
 f_\ald(z)&=&-:\g(z)\left(\hf\g(z)c(z)+C(z)\right)\beta(z):-:c(z)C(z)b(z):\nn
 &+&\kvt\left(\hf\g(z)c(z)+C(z)\right)\pa\var_1(z)-\kvt\left(\hf\g(z)c(z)-
  C(z)\right)\pa\var_2(z)\nn
 &+&\hf(k-1)\pa\g(z)c(z)-\hf(k+1)\pa c(z)\g(z)+k\pa C(z)
\label{Wakbos}
\eea
where we have introduced the simplifying notation
\bea
 \beta(z)&=&\beta_{\al_1}(z)\spa\g(z)=\g^{\al_1}(z)\nn
 b(z)&=&b_{\ald_2}(z)\hspace{.1cm}\spa c(z)=c^{\ald_2}(z)\nn
 B(z)&=&b_\ald(z)\hspace{.2cm}\spa C(z)=c^\ald(z) 
\eea
The energy-momentum tensor is
\bea
 T(z)&=&:\pa\g(z)\beta(z):+:\pa c(z)b(z):+\pa C(z)B(z):\nn
 &+&\hf:\pa\var(z)\cdot\pa\var(z):+\frac{1}{\kvt}\ald_2\cdot\pa^2\var(z)
\eea 
The screening currents of the first kind are 
\bea
 s_1(z)&=&-\left(\beta(z)+\hf c(z)B(z)\right):e^{-\var_1(z)/\kvt}:\nn
 s_2(z)&=&-\left(b(z)-\hf\g(z)B(z)\right):e^{-\var_2(z)/\kvt}:
\eea
{\bf Proposition 3}\\[.2cm]
The screening current $\tilde{s}_1(w)$ is given by 
\ben
 \tilde{s}_1(w)=\left(\beta(w)-\frac{t}{2}c(w)B(w)\right)(\beta(w))^{-t-1}
  :e^{\kvt\var_1(w)}:
\een
and it produces
\bea
 R_{-\al_1,\al_1}&=&t\left(\beta-\frac{1}{2}(t+1)cB\right)\beta^{-t-2}\nn
 R_{-\ald,\al_1}&=&tc\beta^{-t-1}\nn
 R_{-\ald_2,\al_1}&=&R_{\al_1,\al_1}=R_{\ald_2,\al_1}=R_{\ald,\al_1}
  =R_{i,\al_1}=0\nn
 \D(\tilde{s}_1)&=&1
\label{totderosp}
\eea
{\bf Proof}\\[.2cm]
Also in this case the proof is straightforward. One computes all possible 
contractions and obtains after some simple rewritings the total derivatives
(\ref{defR}) and (\ref{totderosp}).\\
$\Box$
 
\section{Conclusion}

In this paper we have continued the search for screening currents of the
second kind in affine current (super-)algebras, and by the explicit 
construction in the case of $SO(5)$ it has been demonstrated that more 
complicated structures than previously anticipated \cite{Ito1,Ras1,PRY4} 
appear in the general case. We intend to come back elsewhere with more 
discussions on the nature of screening currents of the second kind.\\[.5 cm] 
{\bf Acknowledgement}\\[.2cm]
The author is financially supported by the Danish Natural Science Research
Council, contract no. 9700517.

\end{document}